\documentclass{article}

\newcommand{\ds}{\displaystyle}
\def\R{\mathbb R}
\newcommand{\e}{{\rm e}}

\usepackage{amsfonts,amsmath,amsthm,amsxtra,amscd,amssymb,epsfig,graphicx}
\usepackage{color}
\usepackage{chngcntr}
\usepackage{wrapfig}
\usepackage{caption}
\usepackage{subcaption}

\begin{document}
%
%
%
%
%
%
%
%
%

\begin{center}
{\Large Exact nonclassical symmetry solutions of Arrhenius reaction-diffusion.}\\
\end{center}

{\begin{flushleft} P. Broadbridge$^2$, B.H. Bradshaw-Hajek$^1$ and D. Triadis $^{2,3}$.\\
\smallskip

1.  (author to whom correspondence should be sent), School of Information Technology and Math. Sciences, University of South Australia.  Bronwyn.Hajek@unisa.edu.au\\

2. Dept. of Mathematics and Statistics, La Trobe University, Victoria, Australia.\\

3. Institute of Mathematics for Industry, Kyushu University, Japan.
\end{flushleft}

\begin{abstract}
Exact solutions for nonlinear Arrhenius reaction-diffusion are constructed in $n$ dimensions. A single relationship between nonlinear diffusivity and the nonlinear reaction term leads to a nonclassical Lie symmetry whose invariant solutions have a heat flux that is exponential in time (either growth or decay), and satisfying a linear Helmholtz equation in space. This construction extends also to heterogeneous diffusion wherein the nonlinear diffusivity factorises to the product of a function of temperature and a function of position. Example solutions are given with applications to heat conduction in conjunction with either exothermic or endothermic reactions, and to soil-water flow in conjunction with water extraction by a web of plant roots.
\end{abstract}
{\it Keywords:}
Arrhenius, reaction-diffusion, plant root extraction, heat conduction, exact solutions, nonclassical symmetries

\section{Introduction}
Assuming that volumetric heat capacity is constant and neglecting reagent consumption, the temperature  of a reactive mixture satisfies a nonlinear reaction-diffusion equation
\begin{equation}\label{RD}
\ds\frac {\partial \theta ({\bf r},t)}{\partial t}=\nabla . \left[D(\theta)\nabla\theta\right]+R(\theta),\quad
{\bf r}\in \Omega \subset \R ^3,~~t\in[0,t_2),
\end{equation}
where $\Omega$ is compact and connected, $\nabla$ is the usual gradient operator, and $t_2>0.$ $R(\theta)$ is  real-valued at any  value of absolute temperature $\theta\in [0,\infty)$. In order that heat conduction contributes to increasing entropy, $D$ must be positive (e.g. \cite{PBEntropy}). One of the most important forms for the reaction function is the Arrhenius reaction term 
\begin{equation}\label{Arrhenius}
R=\left\{
\begin{array}{ll}
R_0e^{-B/\theta} & \theta>0,\\
0 &\theta=0,
\end{array}\right.
\end{equation}
with $B$ a positive constant. This follows from the Boltzmann-Gibbs equilibrium canonical distribution governing the probability of a particle  overcoming the activation energy barrier $E$ of a reaction. The parameter $B$ may then be identified with $E/k_B$ where $k_B$ is Boltzmann's constant. The rate factor $\rho C R_0$,  $\rho$ being density and $C$ being specific heat at constant volume,  is the amount of heat energy released per unit volume per unit time at very high temperature. This is usually described by a weak power-law dependence on temperature, $R_0=S_0 \theta^m$. For example, kinetic theory of a hard-sphere gas predicts $m=0.5$ \cite{Glassman}. Since the temperature-dependence is dominated by the Arrhenius factor, most extant models assume $m=0$. This will be assumed here, but a straightforward generalisation of the following analysis  can cover the case $m\ne0$ with $\rho C$ also depending on $\theta$.

For first-order reactions, $E$ is  the energy of bond dissociation into reactive molecular components such as free radicals. For example, for the exothermic decomposition of diethyl peroxide, measured reaction rates agree well with the Arrhenius law over a broad range of temperatures \cite{Fine}. In particular, as a combustible mixture is controlled  by extracting heat through the boundaries of the container,  the Arrhenius reaction form is  expected to remain appropriate. 

A full Lie point symmetry classification of (\ref{RD}) was made by Dorodnitsyn {\it et al.} \cite{Dorod}. Various combinations of power-laws and exponential functions for $D(\theta)$ and $R(\theta)$  lead to an expansion of the Lie group of point symmetry transformations of (\ref{RD}) beyond the common Euclidean isometries in space, and translations in time. These in turn open possibilities of invariant solutions that may be obtained by reduction of variables. The invariant solutions indicate a wide range of possible dynamical behaviours, including stable similarity forms for the temperature with multi-peaks, extinguishing, and blow-up in finite or infinite time that can occur even in solutions with compact support \cite{Dorod}.

The non-analytic expression $\exp(-B/\theta)$ is usually approximated by an exponential function, in the Frank-Kamenetskii approximation \cite{FK}, before it can lead to useful exact solutions of (\ref{RD}). Although unbounded reaction terms lead to insight on critical parameters for ignition, there is a significant difference in behaviour when the reaction term is bounded, as in the standard Arrhenius model. The relationship between the full Arrhenius model and the Frank-Kamenetskii model can be more conveniently analysed by using the identity
$$e^{-E/k_B\theta}=e^{-E/k_B\theta_0}\exp\left(\frac{\Theta}{1+\epsilon\Theta}\right);  ~~\epsilon=\frac{k_B\theta_0}{E},$$
where $$\Theta=({\theta-\theta_0}){\frac{E}{k_B\theta_0^2}}$$ is a rescaled temperature with the zero point shifted to some local value $\theta_0$  (e.g. \cite{Fowler}). Wake and Bazley \cite{Wake} extended the lower bound 
$$e^{-E/k_B\theta}\ge e^{-E/k_B\theta_0}\exp\left(\frac{\Theta}{1+\epsilon\Theta_{max}}\right),$$
beyond the maximum temperature to obtain close approximate critical values of parameters at ignition. Gustafson \cite{Gustafson} constructed steady linear and radial solutions of the full nonlinear equation by ``a shooting method combined with a Newtonian-Raphson technique and certain boundary value expansions".

Although the classical Lie point symmetry classification of (\ref{RD}), with $D(\theta)$ and $R(\theta)$ arbitrary, selects only constant or unbounded reaction functions, the complete {\em nonclassical} symmetry classification of reaction-diffusion equations, admits a broader range of possibilities \cite{Arrigo1, Clarkson, Arrigo2, Goard}, even admitting the Arrhenius reaction term \cite{Goard}. Nonclassical symmetries, in the sense of Bluman and Cole \cite{BC}, leave invariant a system consisting of the original governing equation (\ref{RD}) plus the invariant surface condition that restricts the solution set to only invariant solutions. The concept of nonclassical symmetry extends more generally to that of compatibility with an invariance condition, that may be of higher order than a point transformation or a contact transformation \cite{Olver}. Nonclassical symmetry analysis reveals that with a suitable nonlinear diffusivity function, after a change of variables, the Arrhenius reaction-diffusion equation admits separation of variables, resulting in a {\em linear} system. In the following Section 2,  time-dependent radial solutions of this system are constructed in two and three dimensions. In the first example, the complicated nonlinear diffusivity function, as well as the temperature, are given exactly. With the same strongly increasing diffusivity, the exponentially heating solution is mathematically equivalent to that involving conduction through a finite domain with an endothermic reaction and exponential cooling. The strongly increasing diffusivity resembles soil-water diffusivity, with the sink term representing plant root absorption.   In the case of ideal maximal cooling at the boundary, the nonlinear diffusivity is bounded and it is the fixed point of a rapidly converging contraction map. The similarity form of the Kirchhoff variable is given exactly and it is asymptotically of the same form as the temperature. The diffusivity varies so little that the exact solution for the Kirchhoff  variable  closely approximates the temperature at all times.  Particular attention is paid to the case of Newton cooling at the boundary, with non-zero Biot number.

In Section 3, it is shown that the same type of solution construction applies to a class of nonlinear reaction-diffusion equations with spatially varying diffusivity. The results of Sections 2 are extended to allow for spatially variable diffusivity and a nonlinear sink term that applies to water transport in unsaturated soil with extraction by plant roots.

\section{Nonclassical reduction of nonlinear reaction-diffusion to the Helmholtz equation}

A full nonclassical symmetry classification of nonlinear diffusion-reaction equations in two spatial dimensions, was given in \cite{Goard}. Firstly, (\ref{RD}) is expressed in terms of the Kirchhoff variable (e.g. \cite{Philip}),
\begin{equation}
u=u_0+\int_0^{\theta}D(\bar\theta)d\bar\theta~.\label{Kirchhoff}
\end{equation}
If $u_0=-\int_0^{\theta_0}D(\bar\theta)d\bar\theta$ for some $\theta_0 \ge0$, then
$$u=\int_{\theta_0}^{\theta}D(\bar\theta)d\bar\theta~.$$
In that case, a boundary condition $u=0$ corresponds to $\theta=\theta_0$.

In terms of the Kirchhoff variable, the reaction-diffusion equation is
\begin{equation}\label{RDu}
F(u)\frac{\partial u}{\partial t}=\nabla^2u+Q(u),
\end{equation}
where $Q(u)=R(\theta)$ and $F(u)=1/D(\theta)$.
The starting point of the reduction to the Helmholtz equation is the observation that Equation (\ref{RDu}) has a simple nonclassical symmetry 
\begin{equation}\label{nonclass}
\bar u=e^{A\varepsilon}u~;~~\bar t=t+\varepsilon~;~~\bar x^j=x^j
\end{equation}
with $\varepsilon\in \R$, whenever $F$ and $Q$ are related by 
\begin{equation}
Q(u)=AuF(u)+\kappa u, \label{related}
\end{equation}
for some $A,\kappa\in\R$. This is not a classical symmetry because it does not in general leave the equation (\ref{RDu}) invariant, except when one also assumes the invariant surface condition $u_t=Au$. Then the symmetry reduction leads to
\begin{equation}
\nabla^2\Phi+\kappa \Phi=0 \quad {\rm with} \quad u=e^{At}\Phi({\bf x}).
\label{Helmholtz}
\end{equation}

In terms of the original  temperature variable $\theta$, the relation (\ref{related}) is 
\begin{equation}\label{RfromD}
R(\theta)=\left[\kappa+\frac{A}{D(\theta)}\right]\left[u_0+\int_0^\theta D(\bar\theta)~d\bar\theta\right].
\end{equation}
This gives an explicit construction of $R(\theta)$ from $D(\theta)$. Some basic combinations $\left(D(\theta),R(\theta)\right)$, with $R(0)=0$, are given in Table 1. Note that the $R(\theta)$ function in case (d) of Table 1 agrees asymptotically with the Arrhenius reaction term $R_0 e^{-B/\theta}$ as inverse temperature approaches $\infty$.
\begin{table}[h]
\centering
\begin{tabular}{|c|c|c|}
\hline
{}&$\vphantom{\Big(} D(\theta)$ & $R(\theta)$ \rule{0mm}{5mm}\\
\hline
(a)& $\theta^m~(m>-1)$ & $\ds\frac{\kappa}{m+1}\theta^{m+1}+\frac{A}{m+1}\theta$ \rule{0mm}{6mm}\\ 
(b)&$e^\theta$ & $\kappa[e^\theta-1]-A[e^{-\theta}-1]$\rule{0mm}{6mm}\\
(c)&$\cosh(\theta)$ & $\kappa \sinh(\theta)-A\tanh(\theta)$\rule{0mm}{6mm}\\
(d)& $\displaystyle \vphantom{\frac {1}{\Big)}}  \frac {R_0}{\kappa B} \left( 1 + \frac B \theta \right) e^{-B/\theta} - \frac {A}{\kappa} $& 

$\displaystyle \frac {R_0 \theta (B + \theta) e^{-B/\theta}( R_0 e^{-B/\theta}-AB)}{R_0 B (B + \theta) e^{-B/\theta}-AB^2\theta} $
\rule{0mm}{8mm}
\\
\hline

\end{tabular}
\caption{Reaction-diffusion combinations that admit nonclassical scaling symmetry \eqref{nonclass}.}
\end{table}

However, in order to construct $D(\theta)$ from $R(\theta)$, it is necessary to solve a first-order nonlinear differential equation. From (\ref{related}),
\begin{equation}
D(\theta)=u' (\theta)=\frac{Au}{R(\theta)-\kappa u}.
\label{ude}
\end{equation}
By making $u(\theta)$ the subject, then differentiating, this implies
\begin{equation}
AD'(\theta)=D^3\frac{\kappa^2}{R(\theta)}+D^2\kappa~\frac{2A-R'(\theta)}{R}+DA\frac{A-R'(\theta)}{R}.
\label{Dde}
\end{equation}
The  equation \eqref{Helmholtz} for $\Phi$  is simply the linear Helmholtz equation if $\kappa=K^2>0$, and the Laplace equation if $\kappa=0$. 
If $\kappa=-K^2<0$ with $R>0$, from (\ref{ude}) $A$ must be positive. Then the equation for $\Phi$ is the modified Helmholtz equation for which the radial solutions are unbounded either at the origin $r=0$ or at infinity.

\subsection{An explicit solution with Arrhenius reaction: $\kappa=0$}

In the case $\kappa=0$, equation (\ref{ude}) is linear homogeneous, allowing direct integration to obtain
\begin{equation}
u(\theta)=\frac{c_1}{A}\left[\exp\left(\int^{\theta}\frac{A}{R(\bar\theta)}d\bar\theta\right)-1\right],\label{umap}
\end{equation}
with
\begin{equation}
D(\theta)=\frac{c_1}{R(\theta)}\exp\left(\int^\theta \frac{A}{R(\bar\theta)}d\bar\theta\right)\label{Dmap},
\end{equation}
where $c_1$ is an arbitrary  constant that is positive (negative) if $R$ is  a strictly positive (negative) source (sink) term. However the spatial dependence of $u$, given in (\ref{Helmholtz}), is then governed by Laplace's equation. For example, the isotropic radial solutions $u({\bf x},t)=u(r,t)$ in two or three dimensions can only be a constant added to the singular point-source(sink) solutions with the source(sink) strength varying exponentially in time. The isotherms are specified exactly by the mapping that follows from Equation(\ref{umap}):
$$(t,\theta)\mapsto u\mapsto \Phi=e^{-At}u\mapsto r.$$ In two and three dimensions, the radial solutions take the form $r=\exp(\frac{\Phi-c_2}{c_3})$ and $r=\frac{c_3}{\Phi-c_2}$ respectively.\\

For the case of the Arrhenius reaction (\ref{Arrhenius}), the solution (\ref{Helmholtz}) is valid when the diffusivity is exactly
\begin{equation*}
D=\frac{c_1}{R_0}\exp({B/\theta})\exp\left(\frac{A}{R_0}\theta \exp({B/\theta})-\frac{AB}{R_0}E_i(B/\theta)\right). 
\end{equation*}
where $E_i$ is the exponential integral. When $\theta$ is large, we can use the power series for the exponential integral
$$
E_i(x)=\gamma+\ln(x)+\sum_{n=1}^{\infty}\frac{x^n}{nn!}~,
$$
where $\gamma$ = Euler's constant, to find
$$
D\sim\frac{c_1}{R_0}e^{AB(1-\gamma)/R_0}B^{-AB/R_0}\theta^{AB/R_0}e^{A\theta/R_0}~.
$$
When $\theta$ is small, we can use the asymptotic expression
$$
E_i(x)\sim \frac{e^x}{x}\left[1+\sum_{n=1}^{\infty}\frac{n!}{x^n}\right],
$$
to find that $D\to 0$ as $ \theta\to 0.$

This nonlinear diffusivity is plotted in Figure 1b alongside the Arrhenius reaction function, plotted in Figure 1a.
%
%
%
%
\begin{figure}[h]
\centerline{\includegraphics[width=50mm]{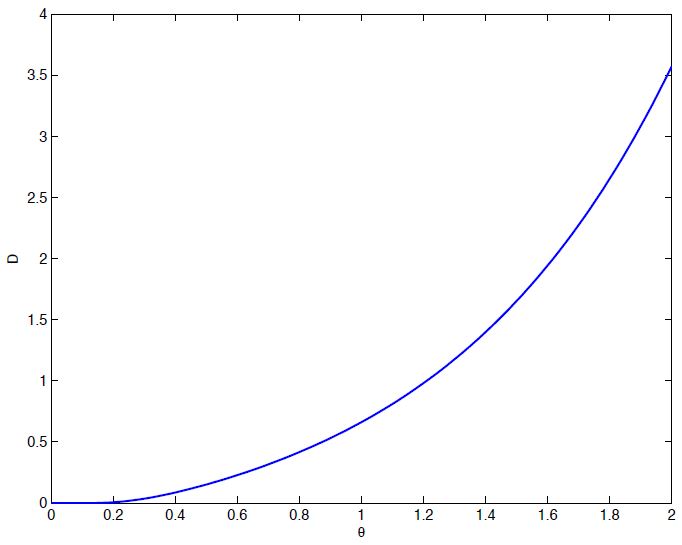}\hspace*{4mm}\includegraphics[width=50mm]{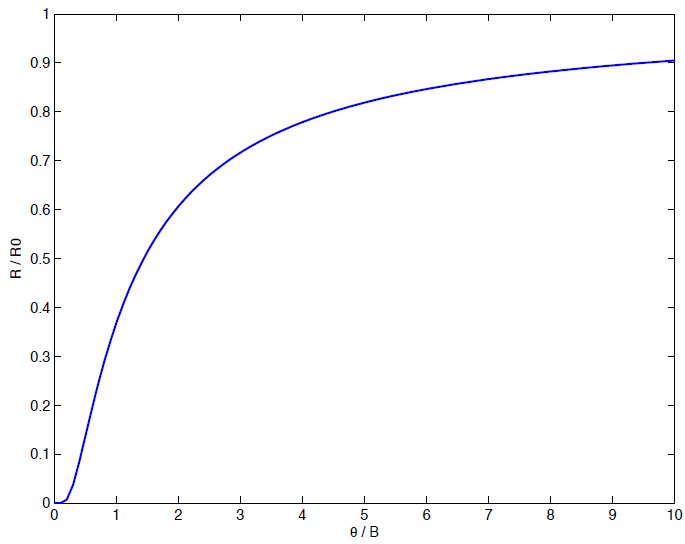}}
\caption{(a) Exact temperature-dependent diffusivity allowing separation of variables of Arrhenius reaction-diffusion. (b) Arrhenius reaction rate, as a function of temperature.}
\label{fig1}
\end{figure}

\noindent{\bf Example 1. A heat conduction problem with  Arrhenius endothermic reaction term}

Even though the enthalpy of reaction is negative, an endothermic reaction will proceed naturally if the reaction products result in a lowering of the Gibbs free energy. An endothermic reaction may have a single activation energy, so that the rate of heat absorption is described by an Arrhenius function. Consider a spherical vessel containing the reagents of an endothermic reaction. $R_0$ is negative.  A small spherical decaying radioactive heat source at the centre, supplies a quantity of heat Q. Then, since $u>0$ and $D> 0$ when $\theta>0$, (\ref{umap})--(\ref{Dmap}) now imply $c_1<0$ and $A<0$. A solution 
\[
u=e^{-|A|t}[c_2-\frac{c_3}{ r}]
\]
can be made to satisfy the boundary conditions
\begin{eqnarray*}
u=0&~{\rm at}~&r=r_1,\\
-4\pi r^2u_r=|A|Qe^{-|A|t} &~{\rm at}~&r=r_0.
\end{eqnarray*}
Note that the latter is simply a heat flux condition that can always be expressed as a linear condition in terms of the Kirchhoff variable $u$. 
\medskip

\noindent{\bf Example 2. Porous media flow with plant-root sink term}

The combination of exponentially growing nonlinear diffusivity, with bounded sink term,  applies to water transport through a web of plant roots \cite{Feddes}. The domain of the problem is considered to be a cylindrical pile of soil of radius $r=r_1$ with a vertical cylindrical injection well in the centre ($r=r_0<r_1$). The dependent variable $\theta$ now designates the water content above the plants' wilting point,  where the sink term approaches zero as roots fail to draw water. In this application, $R_0$ is negative.  Then, since $u>0$ and $D> 0$ when $\theta>0$, (\ref{umap})--(\ref{Dmap}) now imply $c_1<0$ and $A<0$. A solution 
\[
u=e^{-|A|t}[c_2-c_3\ln r]
\]
can be made to satisfy the boundary conditions
\begin{equation}\label{PlantRootbc}
\begin{array}{rclr}
u=0~&{\rm at}&~r=r_1,\\
-2\pi ru_r=|A|Qe^{-|A|t} ~&{\rm at}&~r=r_0.
\end{array}
\end{equation}
These boundary conditions allow for injection of total water volume $Q$ per unit length of injection well into a large cylindrical soil mound that is exposed to the soil-controlled second stage of evaporation (e.g. \cite{Stewart}) at its outer boundary. 

This solution applies analogously to a vertical current-carrying wire along the axis of a cylindrical region, supplying energy for an endothermic reaction.


\subsection{Construction of bounded temperature solutions: $\kappa>0$}

\noindent{\bf Example 3. Heat conduction with exothermic reaction in a compact region}

With $\kappa=K^2>0$, $u=\Phi({\bf r})\exp(At)$ satisfies the linear Helmholtz equation. A non-negative isotropic solution can satisfy boundary conditions
\begin{eqnarray*}
u_r=0,~~r=0 , \\
u=0,~~r=r_1,
\end{eqnarray*}
by choosing $\Phi=j_0(Kr)$ in three dimensions, $\Phi=J_0(Kr)$ in two dimensions and $\Phi=\cos(Kr)$ in one spatial dimension, with $K=\lambda_1/r_1$, $\lambda_1$ being the first zero of the spherical Bessel function $j_0$  ($\lambda_1=\pi$), the standard Bessel function $J_0$ ($\lambda_1\approx 2.4048$) or the cosine function ($\lambda_1=\pi/2$) respectively. Further, from the separation of variables in (\ref{Helmholtz}), that same solution satisfies other linear homogeneous boundary conditions such as
$$-u_r=Bi~u,~~~r=r_2<r_1,$$ with $Bi$ constant. This resembles the linear condition for Newton cooling, in which $Bi$ is the Biot number, after rescaling and non-dimensionalising variables. The left hand side is heat flux but the right hand side is a multiple of the Kirchhoff variable, rather than temperature. Below, it is demonstrated that the nonlinear diffusivity has bounded variation and is close to being constant. This means that the above boundary condition is very close to that of Newtonian cooling to a very low-temperature environment. If the physical parameters $r_2, D(0)$ and $B_i$ are prescribed, then the solution parameters are given by $A=-K^2D(0)$,
where $K$ is the unique solution of $-\Phi'(Kr_2)/\Phi=B_i/K$. $K$ must lie between 0, where $ -\Phi'(Kr_2)/\Phi=0$ and $\lambda_1/r_2$, where $-\Phi'(Kr_2)/\Phi$ approaches infinity. With these homogeneous boundary conditions, there always exists an extinguishing solution in which the temperature approaches zero uniformly everywhere. A demonstration of the stability of this similarity solution is given in the Appendix A. \\

In the case $\kappa>0$, (\ref{ude}) is equivalent to the canonical form for an Abel equation of the second kind, via
\begin{eqnarray}
w=\kappa u-R(\theta),~~z=-A\theta -R(\theta),\nonumber\\
w\frac{dw}{dz}=w+\Phi(z) ;~~\Phi(z)=\frac{AR(\theta)}{A+R'(\theta)}.
\label{Abel2}
\end{eqnarray}
The standard list of known integrable forms of the Abel equation is given in \cite{Polyanin}. In principle, this list may be used to furnish $R(\theta)$ functions that produce solvable forms of equation (\ref{Abel2}). If $g(z)$ a particular solution for an integrable form of (\ref{Abel2}), and $g^{-1}(w)$ is the corresponding inverse function, an $R(\theta)$ function that produces (\ref{Abel2}) for $w(z)$ is given implicitly according to 
\begin{equation}
\theta = \frac {R + g^{-1}(-R)}{-A}.
\end{equation}
 It is not clear if the segmented solution method of  \cite{Panay} is practicable for the general case.\\

A good approximate analytic reconstruction of $D(\theta)$ may be obtained by applying a contraction map towards a fixed point.
With $u_0=0$, the solution to (\ref{ude}) satisfying initial value $D(\theta)\to0,~~\theta\to 0,$ must be a fixed point of the  map
$$
D_{n+1}(\theta)={\mathcal M}D_n(\theta)=\frac{-A\int_0^{\theta}D_n(s)ds}{K^2\int_0^{\theta}D_n(s)ds -R(\theta)}.$$
When considering the separated solution (\ref{Helmholtz}), $|A|$ and $K$ may be set to $1$ by choosing dimensionless length and time coordinates $Kr$ and $|A|t$.
\begin{equation}
D_{n+1}(\theta)=\frac{\bar D_n(\theta)}{\bar D_n(\theta) -R(\theta)/\theta},
\label{iterate}
\end{equation}
where $\bar D_n(\theta)$ is the running mean value, 
\begin{equation}\label{mean}
\bar D_n(\theta)=\frac{1}{\theta}\int_0^{\theta}D_n(s)ds.
\end{equation}

In the case of the Arrhenius reaction term, we may set $R(\theta)$ to be $R_0e^{-1/\theta}$ by using a dimensionless temperature variable $\theta/B$. Then the maximum value of $R(\theta)/\theta$ is $R_0/e.$ \\

From (\ref{iterate}), using $\bar D(0)=D(0)$, it follows that $D(0)=1$. From (\ref{iterate}) it also follows that $D(\theta)> 1$ for $\theta>0$. If $D(\theta)$ is bounded for $\theta\in(0,\infty)$ and monotonic for $\theta$ sufficiently large, $D(\theta)$ must have a limit $D(\infty)=\bar D(\infty)$. In that case it follows  from (\ref{iterate}) that $D(\infty)=1$.\\

Consider $\mathcal M$ as a mapping on the set ${\mathcal S}$ of bounded continuous functions $f$ on $(0,\infty)$ such that $f-1$ is non-negative,
$${\mathcal S}=\{f\in C^1(0,\infty):\exists b~\forall x\in (0,\infty)~ 0\le f(x)-1\le b\}.$$ Suppose that $D_0,E_0\in\{{\mathcal S}\}$, $D_n=M^nD_0,~E_n=M^nE_0$, where $E_{n+1}(\theta)$ and $\bar E_n(\theta)$ are defined in a analogous way to \eqref{iterate} and \eqref{mean} above. Now 
\begin{eqnarray*}
\left| D_{n+1}(\theta)-E_{n+1}(\theta)\right|&=&\left|\frac{\int_0^{\theta} D_n(x)dx}{\int_0^{\theta} D_n(x)dx-R(\theta)}-
\frac{\int_0^{\theta} E_n(x)dx}{\int _0^{\theta}E_n(x)dx-R(\theta)} \right|\\
&=&\left| \frac{\theta{\bar D_n}}{\rule{0mm}{3.5mm}\theta{\bar D_n}-R(\theta)} -\frac{\theta{\bar E_n}}{\rule{0mm}{3.5mm}\theta{\bar E_n}-R(\theta)}\right|\\
&=&\left| \frac{({\bar D_n}-{\bar E_n})R(\theta)/\theta}{\rule{0mm}{3.5mm}[\bar D_n-R(\theta)/\theta][\bar E_n-R(\theta)/\theta]}\right|.
\end{eqnarray*}
 
If $R_0<e$, then $\inf\{D(\theta)\}=1>\sup\{R(\theta)/\theta\}=R_0/e$, consequently
$$
\frac{\|E_{n+1}-D_{n+1}\|_{\infty}}{\|E_n-D_n\|_{\infty}}\le\frac{\sup\{R(\theta)/\theta\}}{[\inf\{D_n(\theta)\}-\sup\{R(\theta)/\theta\}][\inf\{E_n(\theta)\}-\sup\{R(\theta)/\theta\}]}.
$$
\begin{equation*}
\Rightarrow \frac{\|E_{n+1}-D_{n+1}\|_{\infty}}{\|E_n-D_n\|_{\infty}}\le\frac{1}{e/R_0-2+R_0/e}.
\end{equation*}
This is a contraction map, provided $$R_0<\frac{(3-\sqrt 5)e}{2}\approx 1.0383~.$$
For example, consider the case $R_0=1$ for which the above estimate of the contraction factor is $1/1.086$. In practice, the convergence rate is better than this modest value suggests.
In the case $R(\theta)=\exp(-1/\theta)$, since $D(0)=D(\infty)=1$, it is natural to choose $D_0(\theta)=1$.  Using the standard notation $\beta=1/\theta$, the iterative estimates are
\begin{eqnarray}
\nonumber
D_0&=&1,\\
\nonumber
D_1&=&\frac{1}{1-\beta \e^{-\beta}}=\sum_{n=0}^{\infty}\beta^ne^{-n\beta}~
\hbox{(noting $\beta e^{-\beta}\le e^{-1}<1$)},\\
\nonumber
D_2&=&\frac{I}{I-\e^{-\beta}};\\
\nonumber
I&=&\int_\beta ^{\infty} \frac {D_1}{\beta^2} d\beta=\beta^{-1}-E_i(-\beta)+
\sum_{n=2}^\infty \frac {(n-2)!}{n^{n-1}} e^{-n\beta} \sum_{k=0}^{n-2} \frac {n^k \beta^k}{k!} 
\end{eqnarray}
Each term of the $I$ summation is bounded above by ${(n-2)!}/{n^{n-1}}$, which is of order $e^{-n}/\sqrt n$ as $n \to \infty$. The $n=5$, $e^{-5\beta}$ terms have a combined value of less than 0.01.

Note that the sum of the $e^{-5\beta}$ terms has a maximum value of less than 0.01.\\
In Appendix B, $D(\theta)$ is calculated accurately after obtaining exact series forms for $u(\theta)$. The simpler approximate diffusivity functions $D_j(\theta)$ are shown in Figure \ref{fig:2} to agree well with the exact solution. 
\begin{figure}
\includegraphics[width=120mm]{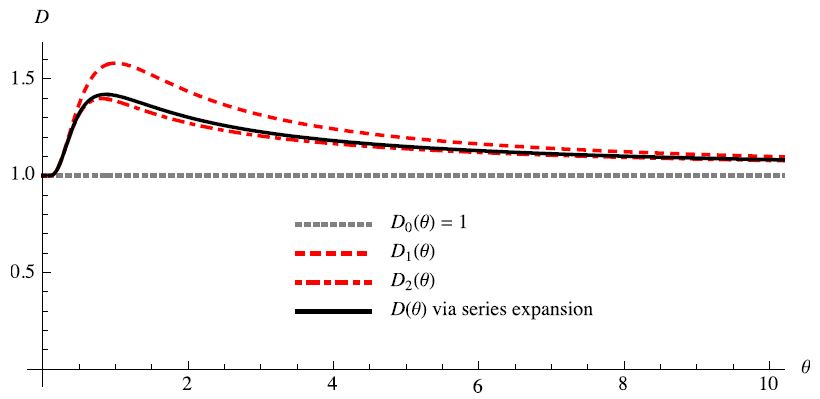}
\caption{$D(\theta)$ constructed by series about $\theta=0,0.5,10$, as well as approximations $D_1(\theta), ~D_2(\theta)$ }
\label{fig:2} 
\end{figure}

The diffusivity function has a single local maximum where its value is no higher than 45\% above its mean value of 1.  { From (\ref{Dde}), a stationary point of $D(\theta)$ may occur only where $D(\theta)=1+R'(\theta)=1+R_0\theta^{-2}e^{-1/\theta}$. This expression has a maximum value, implying 
\begin{equation}D_m\le 1+4R_0/e^2\approx 1.5413~.
\label{estDm}
\end{equation}
}\\

In effect, since at small-$t$, local disturbances extend by diffusion to a depth proportional to $\sqrt{Dt}$ and $\sqrt D$ varies by only 20\%, the diffusivity is effectively constant for some practical purposes. This is seen in the solution, plotted in Figure 3, wherein the temperature $\theta(r,t)$ asymptotically approaches the Kirchhoff variable $u(r,t)$ which is identical to the temperature when $D=1$.\\
\begin{figure}
\centerline{\includegraphics[width=100mm]{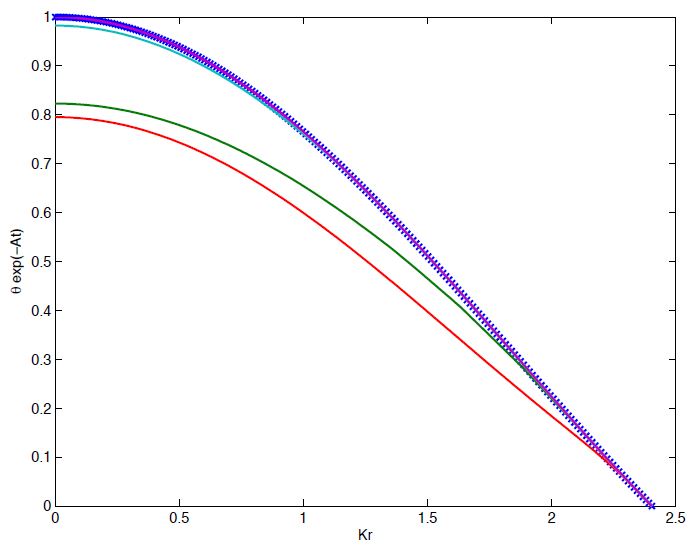}}
\caption{Temperature profiles at times $|A|t$=-1.5,0,1.5 and 2.5 (solid) approaching the Bessel function profile of u(r,t) (x crosses) from below, as time increases.}
\label{fig:3} 
\end{figure}

The construction of the nonlinear diffusivity function by a contraction map relied on the fact that with Arrhenius reaction,  $\beta Q(\beta)~(=R(\theta)/\theta)$ has an upper bound less than $1$. The same construction will apply to other reaction laws with this property, for example $Q(\beta)=\beta^m\exp(-\beta)$ with $m>-1$. This includes the case $m=-1/2$ that follows from the kinetic theory of gases.

\subsection{The case $\kappa<0$.}

If $\kappa=-K^2<0$, from (\ref{ude}) with $R(\theta)>0, ~A$ must be positive. The diffusivity, $D(\theta)$, has a singularity where at leading order $D'(\theta)\sim(\kappa^2/A)D^3$ and $D\sim (\theta_a-\theta)^{-1/2}$
at  some freely chosen positive  value $\theta =\theta_a$. In the modelling of unsaturated soil-water flow, the location of the singularity is chosen to be slightly greater than the volumetric water content at saturation \cite{BW2}.

Also the equation for $\Phi$ is the modified Helmholtz equation for which the radial solutions are unbounded either at the origin $r=0$ or at infinity. There is an exact solution for exothermic reaction-diffusion on the domain $r>r_0$, 
\begin{equation}
u=e^{At}K_0(Kr),
\label{modBess}
\end{equation}
where $K_0$ is the modified Bessel function of the second kind. Although $u$ is unbounded, the corresponding value of $\theta$ remains less than the singular value $\theta_a$, due to the singularity in $D(\theta)$ within the integrand in the Kirchhoff transformation $\theta\to u$.  This singular nonlinear diffusivity does not occur in heat conduction but it may occur in models of unsaturated soil-water flow. However, in that application, positive distributed water sources are not common. It is more common to consider distributed sinks, due to plant roots.\\

\noindent{\bf Example 4. Porous media flow with plant-root sink term: $\kappa<0$}

As in Example 2, we consider a distributed sink with $R_0<0$, due to plant roots. Because of plant root extraction, $A<0$ so that a finite amount of water is supplied over all time through the inner surface at $r=r_0$. From (\ref{RfromD}), $D>0, A<0$ and $\kappa<0$ implies $R<0$. However in this case, from the solution (\ref{modBess}), $\theta$ does not reach zero at any finite value of $r$ but instead approaches zero at infinite distance. It is to be noted that the modified Helmholtz equation has previously been applied to {\em steady state} solutions of unsaturated soil-water flow, by analogy to external problems of Helmholtz acoustic scattering \cite{Waechter}. Just as in the latter, we may construct analogous exterior solutions for boundaries of various shapes. However, at large values of $r$ it is well known that the solutions agree asymptotically with the isotropic ones, so the spherical or circular boundaries are canonical.

\section{Incorporating heterogeneity}

{We now consider the case in which the medium (or substrate) is spatially heterogeneous. The appropriate governing reaction-diffusion equation can be written 
\begin{equation}\label{hetero_t}
\theta_t=\nabla.[f({\bf x})D(\theta)\nabla\theta]+R(\theta),
\end{equation}
where ${\bf x}\in\Omega\subset\R^n$.} {Here, the heterogeneity is represented by a positive differentiable function $f$ that is the amplification factor of the diffusivity. By rewriting this equation in terms of the Kirchhoff variable \eqref{Kirchhoff}, we obtain
\[
F(u)u_t=\nabla.[f({\bf x})\nabla u]+Q(u),
\]
(the heterogeneous  equation corresponding to \eqref{RDu}) where $Q(u)=R(\theta)$ and $F(u)=1/D(\theta)$, as before. Equation \eqref{hetero_t} admits the same simple nonclassical symmetry \eqref{nonclass}, whenever $F(u)$ and $Q(u)$ are related by Equation \eqref{related}.{This means that $R(\theta)$ can be constructed from $D(\theta)$ using \eqref{RfromD}, or $D(\theta)$ could be constructed from $R(\theta)$ by solving \eqref{ude}. By again setting $u({\bf x},t)=\e^{At}\Phi({\bf x})$ we derive the second order linear equation for $\Phi(\bf x)$,
\begin{equation}\label{h_phi_ode}
f({\bf x})\nabla^2\Phi({\bf x})+\nabla f({\bf x}).\nabla\Phi({\bf x})+\kappa\Phi({\bf x})=0,
\end{equation}
which is valid in $n$ dimensions. In the case where $A<0$, the solution for $\theta$ will be asymptotically similar to $u$. The relationship between $D(\theta)$ and $R(\theta)$ is the same as that described in the preceeding section, and the solutions for $D(\theta)$ found about also hold in the hetergeneous examples described below.

\subsection{The case $\kappa=0$}

\noindent{\bf Example 5. Porous media flow with a heterogeneous substrate: $\kappa=0$}

As an example, we now reconsider the porous medium plant-root extraction problem described in Section 2.1 above, with the inclusion of spatial heterogeneity representing a changing scale of soil pores. In Miller self-similar soils, $f$ is the geometric scale factor of the soil pores (\cite{Miller,PBJMP}), with $f({\bf r}_0)=1$ at some chosen reference point. 

In this example, the problem is written in (one-dimensional) radially symmetric coordinates and the heterogeneity is described by $f({\bf x})=f(r)$. The differential equation for {$\Phi(r)$ replacing $\Phi({\bf x})$ in \eqref{h_phi_ode}, is}
\[
\ds\frac{f}{r}\frac{d}{d r}\left(r\ds\frac{d\Phi}{dr}\right) +\ds\frac{df}{dr}\ds\frac{d\Phi}{dr}=0,
\]
which can be directly integrated for arbitrary $f(r)$ to give

$$u(r,t)=e^{At}\Phi(r)=ce^{At}\int\frac{1}{rf(r)}dr$$
with $c$ constant.
   From (\ref{PlantRootbc}a), we deduce that $A<0$. Having freely chosen $f(r_0)$ to be $1$, the solution, in terms of the Kirchhoff variable, satisfying  boundary conditions \eqref{PlantRootbc}, is
$$u(r,t)=\frac{|A|}{2\pi}e^{At}\int_r^{r_1}\frac{1}{sf(s)}ds.$$

\subsection{The case $\kappa>0$}

\noindent{\bf Example 6. Heat conduction in a heterogeneous medium}

Let $\kappa=K^2$. The differential equation for $\Phi(r)$ replacing $\Phi({\bf x})$ in \eqref{h_phi_ode} is
\begin{equation}\label{HeteroKpos}
\ds\frac{f}{r}\frac{d}{d r}\left(r\ds\frac{d\Phi}{dr}\right) +\ds\frac{df}{dr}\ds\frac{d\Phi}{dr}+\kappa\Phi=0.
\end{equation}
With $f(r)$ any power-law of the cylindrical radius, exact solutions are readily available (e.g. \cite{Polyanin}).

With heterogeneity described by the factor $f(r)=r_0/r$, \eqref{HeteroKpos} reduces to Airy's equation, with general solution
$$\Phi=c_1A_i(-K^{2/3}r_0^{-1/3}r)+c_2B_i(-K^{2/3}r_0^{-1/3}r).$$  

When $f(r)$ is of the form
\begin{equation}\label{hetero1}
f(r)={\left(\frac{r}{r_0}\right)^2},
\end{equation}
the resulting ODE is merely a homogeneous Euler equation, with general solution
\begin{equation*}
\Phi(r)=c_1r^{-1+\sqrt{1-(K r_0)^2}}+c_2r^{-1-\sqrt{1-(K r_0)^2}}.
\end{equation*}
This solution can be made to satisfy boundary condition $\Phi(r_1)=0$ by setting $c_2=-c_1r_1^{2\sqrt{1-(Kr_0)^2}}$.\\
For the case with $Kr_0>1$, the general real solution is
\begin{equation*}
\Phi(r)=\frac{1}{r}\left[c_1\cos(\omega\log r)+c_2\sin(\omega\log r)\right];~~~\omega=\sqrt{(Kr_0)^2-1}.
\end{equation*}
For the above cases, solution parameters may be chosen so that $u(r_1)=0$. However, $r=0$ is a singular point where thermal conductivity is either zero or infinite. The solution may be regarded as exterior to the surface of a hot wire at $r=r_0$ from where a finite quantity of heat is supplied.}

\subsection{The case $\kappa<0$}

\noindent{\bf Example 7. Porous media flow with a heterogeneous substrate: $\kappa<0$}

With the same heterogeneity as in the previous model \eqref{hetero1} but with $\kappa=-K^2<0$, the solution with $f(r)=(r/r_0)^2$ is simply
\begin{equation*}
\Phi(r)=c_1r^{-1+\sqrt{1+(Kr_0)^2}}+c_2r^{-1-\sqrt{1+(Kr_0)^2}}
\end{equation*}
and with $f(r)=r_0/r$,
$$\Phi=c_1A_i(K^{2/3}r_0^{-1/3}r)+c_2B_i(K^{2/3}r_0^{-1/3}r).$$ 
The parameters $c_i$ and $K$ may be chosen so that the exterior solution satisfies the same boundary conditions on $[r_0,\infty]$ as in the case of a homogeneous medium.

\section{Conclusion}

Provided the nonlinear diffusivity and the nonlinear reaction term satisfy a single relationship, the Kirchhoff variable $u$, which is the integral of nonlinear diffusivity, admits solutions that are obtainable by separation of variables to a linear system, whose solution is an exponential in time, multiplying an arbitrary solution of the Helmholtz, modified Helmholtz or Laplace equation in space. The heat flux is merely $-\nabla u$, which is given explicitly everywhere in these solutions.

 If the nonlinear diffusivity function is specified, then the compatible reaction function can be constructed directly by integration. If the nonlinear reaction function is specified, then the diffusivity function is the solution of a differential equation that is equivalent to an Abel equation if u satisfies a Helmholtz equation, or to a separable equation if $u$ satisfies Laplace's equation. To the best of our knowledge, this provides the only known exact closed-form solutions in two and three dimensions, of nonlinear reaction-diffusion equations with the classical Arrhenius reaction term. Even when $D(\theta)$ satisfies an Abel equation that is not of known solvable type, in some circumstances it may be specified to arbitrary precision using exact series expansions. 
 
This construction is valid in any natural number $n$ of dimensions, and it generalises also to heterogeneous extensions of the Helmholtz factor equation. Applications are given for radial solutions of exothermic reactions with nonlinear heat conduction, endothermic reactions with nonlinear heat conduction, and water flow from a supply well into a cylindrical soil mound with soil-limited evaporation at the outer boundary. The logistic shape of the negative Arrhenius reaction term resembles the behaviour of distributed plant roots that have a maximum and minimum value of water uptake rate near saturation and wilting point respectively.

As is well known from acoustic scattering theory, exterior solutions of the Helmholtz equation  typically asymptotically approach radial symmetric solutions at large distances from the scattering surface. Hence, the radial solutions illustrated here are in a sense, canonical. The solution method used here, involves a free function of the Helmholtz equation so it would not be difficult to use known non-radial solutions. The special solutions presented here could at least be used as bench tests for more broadly applicable approximate solution methods.

This approach will lead to ongoing investigations of other nonlinear partial differential equations with more than one free function, that may admit special nonclassical symmetry reductions.\\

\medskip

{\bf Conflicting interests:} The authors have no conflicting interests.\\

{\bf Authors' Contribution}
PB conceived the general approach, provided the solution in the case K=0, provided recursive approximations for diffusivity  and directed the project. BH provided the solutions for heterogeneous models. DT provided series constructions for nonlinear diffusivity. All authors contributed to the writing, and all agree on the current version.\\

{\bf Acknowlegements.}
We thank Dr Joanna Goard, Prof. Joel Moitsheki  and Dr Edoardo Daly for useful discussions on the symmetry reductions and plant-root absorption.

 \appendix{\bf Appendix A. Stability of similarity solution.\\}
 It is shown here that the extinguishing radial similarity solution of the Arrhenius reaction-diffusion equation given in Section {2.2}, is exponentially stable to small perturbations. Change variables to a set of canonical variables of the nonclassical symmetry, namely $r,t$ and $v=ue^{-At}$. In terms of these variables, the reaction-diffusion equation may be expressed
 $$F(u)v_t=\nabla^2v+K^2v,$$
 with boundary condition $v(r_1)=0.$
 The similarity solution is a pseudo-steady state, $v_s=\Phi(r)=A_0J_0(\lambda_1r/r_1)$, corresponding to exponential decrease of the Kirchhoff variable $u_s=\Phi(r)e^{-|A|t}$ and satisfying $\Phi '(0)=0$ and $ \Phi (r_1)=0.$ Now consider a perturbed solution, in plane polar coordinates $$v=\Phi(r)+w(r,\phi,t),$$ with twice-differentiable initial condition $$v=\Phi(r)+w_0(r,\phi);~ |w_0|<\epsilon <<1.$$ Then order-$\epsilon$ perturbation $w$ satisfies the same homogeneous boundary conditions, plus $2\pi$-periodicity with respect to $\phi$, plus
 $$F\left(u_s(r,t)\right)w_t=\nabla^2w+K^2w. $$ Unlike the original equation for temperature $\theta$, this is a linear equation with no squared derivative terms, allowing recourse to comparison theory of reaction-diffusion equations (e.g. \cite{Ge}). Note that $1/F(u)=D(\theta)$ and that from Section 2.2, $D_m\ge D\ge D(0)=-A/K^2$. Hence by comparison,  $|w(r,t)|$ must decay at least as fast as the solution to the linear initial-boundary problem with smaller diffusivity and larger reaction term,
 \begin{equation}
 Q_t=D(0)\nabla^2Q+D_mK^2Q~~;~~Q(r,\phi,0)=w_0(r,\phi);~~Q(r_1,\phi,t)=0.
 \nonumber
 \end{equation}
 Let 
 $$\tau=|A|t/\kappa^2; \hbox{ and }q=Qe^{-[D_m/D(0)]\kappa^2\tau}.$$
 Then we have the standard linear heat equation $q_\tau=\nabla^2 q$ with $q$ satisfying the same initial and boundary conditions as $Q$.
 The solution $q(r,\phi,t)$ may be expanded as a standard Fourier-Bessel series (e.g. \cite{Ozisik}). Without loss of generality, we neglect the first component in the series expansion of $w_0$, which is a multiple of $J_0(\lambda_1 r/r_1)$, and which takes the form of the assumed unperturbed similarity solution for $v$. Each of the other terms in the series for $Q$ is of 
 the form 
 \begin{equation}
 J_n(\lambda_{n,m}r/r_1)\left[A_{n,m}\cos (n\phi)+B_{n,m}\sin(n\phi)\right]e^{([D_m/D(0)]\lambda_{0,1}^2-\lambda_{m,n}^2)\tau /r_1^2},
 \label{Jn}
 \end{equation}
 where $n\in\mathbb N$, {$\lambda_{n,m}$} is the $m$'th zero of Bessel function $J_n$, while {$A_{n,m}$ and $B_{n,m}$} are arbitrary real coefficients. Note that beyond $\lambda_{0,1}$, the next smallest root is $\lambda_{1,1}$. From the estimate (\ref{estDm}),
 each of the terms (\ref{Jn}) is decreasing exponentially in time, provided
 \begin{equation}
 R_0<e^2[(\lambda_{1,1}/\lambda_{0,1})^2-1]/4\approx 2.8425.
 \end{equation}
 In terms of the original physical parameters of the unscaled boundary value problem, this criterion is
 \begin{equation}
 \frac{R_0r_1^2}{B\lambda_{0,1}^2D(0)}~<2.8425,
\end{equation} 
which is satisfied in practical cases. \\

\appendix{\bf Appendix B. Evaluating $D(\theta)$ using asymptotic series.\\}

 The reaction term $R(\theta) = R_0 \exp( - B/\theta)$ is of special interest.
With $\kappa = K^2$, the parameters $K$, $|A|$ and $B$ in equation (\ref{ude}) may be set to 1 by adopting appropriate dimensionless variables.  By successive isolation of leading-order terms, the asymptotic behaviour of $u(\theta)$ near $\theta = 0$ can be shown to be
\begin{equation}
u(\theta) \sim - {\rm sign}(A) \theta + R_0 \theta \exp(-1/\theta) \sum_{r=0}^\infty \theta^{-r}\exp(-r/\theta)
\sum_{m=0}^{\infty} \theta^m q_{r,m}.
\label{th0assform}
\end{equation}
This series structure assumes $u_0 =0$ but extension to $u_0 > 0$ should be straightforward.
Adopting $q_{r,-1} = 0$, substitution into (\ref{ude}) shows that
\begin{align}
q_{0,m} &= (-1)^m m!,\quad {\rm for} \quad m \ge 0;\\
q_{r,0} &= \frac {(-{\rm sign}(A) R_0)^r}{r+1},\quad {\rm for} \quad r \ge 0;\\
q_{r+1,m+1} &= \frac{ {\rm sign}(A) R_0}{r+2} \bigg\{ \frac {{\rm sign}(A) (r-m)}{R_0} q_{r+1,m} + (r-j-1) q_{r,m} 
\nonumber \\
&\quad - 
(r+1)q_{r,m+1} 
+\sum_{s=0}^r \sum_{l=0}^m q_{r-s,m-l}\big[(s+1) q_{s,l}- (s-l)q_{s,l-1} \big] \bigg\},
\nonumber \\
&\quad {\rm for} \quad r,m \ge 0.
\end{align}
With the set of coefficients $\{q_{r,m}\}$ determined iteratively as above, the series (\ref{th0assform}) is divergent, but can still be used to give accurate specifications of $u(\theta)$ for $\theta$ sufficiently small. Partial sums of (\ref{th0assform}) may be produced by evaluating both $r$ and $m$ sums up to some maximum integer. Performance of the partial sums is  improved by substituting the known asymptotic form
\begin{equation}
- \frac {\exp(1/x)}{x} {\rm Ei}(-1/x) \sim  \sum_{m=0}^\infty m! (-x)^m.
\end{equation}

A Taylor series expansion for $u(\theta)$ centered on a non-singular point \mbox{$0 < \theta_{0} < \infty$} is comparatively easy to derive. We adopt 
\begin{equation}
u(y) = \sum_{n=0}^\infty \lambda_n y^n, \qquad y = \frac { \theta_{0} - \theta}{\theta_{0}};
\end{equation} 
and by substitution, can show that
\begin{align}
\lambda_{i+1} &= \frac 1{(i+1)(\lambda_0 - R_0 \exp(-1/u_{0}))}\bigg\{ {\rm sign}(A) \theta_{0} \lambda_i 
\nonumber \\
& \quad - \sum_{n=1}^i (i-n+1) \lambda_{i-n+1} \Big[\lambda_n + 
\frac {R_0}{\theta_{0}} {}_1F_1\big(  n+1, 2; -1/ \theta_{0} \big) \Big] \bigg\},
\\
&\quad {\rm for} \quad i \ge 0.
\nonumber
\end{align}
Here $\lambda_0 = u(\theta_{0})$ is assumed known, and the above sum evaluates to zero for $i = 0$. 
The asymptotic series form of $u(\theta)$ as $\theta \to \infty$ can also be ascertained and appears to have a non-zero radius of convergence.\\

We can use the above series to evaluate $u$ for $ 0 < \theta < \theta_{max}$ to a prescribed accuracy by matching expansions about different points. To demonstrate, take $R_0 =1$ and $A < 0$. 
Using (\ref{th0assform}) with $r_{max} = m_{max} = 10$, we can ascertain $u(1/10) =  0.1000041579094$ to $12$ significant figures.
A Taylor series expansion about $\theta_{0} = 1/2$, with $u(1/2) = $$0.55582409195937$ then covers the domain $0.1 < \theta < 0.9$, and matches the accuracy of the known value at $\theta = 1/10$. Introducing a second expansion about  $\theta_{0} = 10$, with $u(10) = 11.79313028084656$, extends the domain of the solution to $\theta = 19.5$, and is compatible with the known value $u(9/10) = 1.1135172087801$. Figure 2 shows the resulting diffusivity $D(\theta) = d u/d\theta$. To obtain a solution of greater accuracy, we can begin with evaluation at a point $\theta < 1/10$ according to series (\ref{th0assform}). A greater number of linked Taylor series may then be needed to cover $0 \le \theta_* \le 19.5$ as done above.

\end{document}